\documentclass{article}
\usepackage{amsfonts}
\newtheorem{theorem}{Theorem}
\newtheorem{corollary}{Corollary}

\begin{document}

\author{Iskander A. TAIMANOV
\thanks{Institute of Mathematics, 630090 Novosibirsk, Russia;
taimanov@math.nsc.ru. The author was supported by the Russian Foundation
for Basic Researches and the Science Support Foundation.}
}
\title{Tame integrals of motion and 0-minimal structures}
\date{}
\maketitle

\def\R{\mathbb R}

\section{Introduction}

An integrability of a Hamiltonian flow presumes that we can describe the
flow in a simple way and, in particular, the decomposition of its
phase space into invariant Liouville tori and the singular locus
looks very simple geometrically.  In this case it is said that
integrals of motion are tame. Usually this is the case when the
flow in integrable in terms of (real-)
analytic functions as it was already shown
for different examples. Here we would like
to expose some ideas from mathematical
logics, i.e. from the theory of o-minimal structures,
which clarify an analytic notion of a tame integral of
motion, and demonstrate this approach for
the integrability problem of geodesic flows.

\section{Different meanings of integrability}

Let $M^{2n}$ be a symplectic manifold with a symplectic form
$$
\omega = \sum_{i <  j} \omega_{ij} dx^i \land dx^j.
$$
Let us denote by
$\omega^{ij}$ the inverse matrix to a skew-symmetric matrix
$\omega_{ij}$. Then any smooth function $H$ on this manifold define
a Hamiltonian flow by the equation which describes the evolution of
any smooth function $f$ along the trajectories of the flow:
$$
\frac{df}{dt} = \{f,H\} = \omega^{ij}
\frac{\partial f}{\partial x^i}\frac{\partial H}{\partial x^j}.
$$
We recall that the function $H$ is called
the Hamiltonian function of the flow (or just the Hamiltonian)
and the skew-symmetric operation $\{f,g\}$ on smooth function is called the
Poisson brackets.

It is said that a function $f$ is a first integral, or an integral
of motion, of the flow if it is preserved by the flow:
$$
\{f,H\} = 0.
$$

The Liouville or complete integrability of a Hamiltonian flow is defined
as follows:

\begin{itemize}
\item
a flow on an $n$-dimensional symplectic manifold $M^{2n}$
is called integrable if it admits a family of first integrals
$I_1,\dots,I_n = H$ such that these integrals are in involution:
$$
\{I_j,I_k\} = 0, \ \ 1 \leq j,k \leq n,
$$
and they are functionally independent almost everywhere, i.e. outside
some nowhere dense set $\Sigma$ which is called the singular locus.
\end{itemize}

The mapping
$$
\Phi: M^{2n} \to \R^n, \ \ \Phi(q) = (I_1(q),\dots,I_n(q)),
$$
is called the momentum map.

The condition that the first integrals are in
involution means that Hamiltonian
flows generated by the functions $I_1,\dots,I_n$ as Hamiltonians commute
everywhere. If $X$ is a compact component of the level surface $\Phi = c$
on which the integrals are functionally independent then this component is
diffeomorphic to a torus on which the Hamiltonian flows corresponding to
$I_1,\dots,I_n$ are linearized and moreover such a linearization can be
extended to a neighborhood of $X$.

The proofs of these statements on integrable flows are exposed
in \cite{Arnold}.

But we see that there is some freedom in the definition when we say about
functional independence of first integrals almost everywhere.
It could be that

\begin{itemize}
\item
they are functionally independent on an open dense set;

\item
given a smooth measure on $M^{2n}$
such that the measure of $M^{2n}$ is finite, the first integrals are
functionally independent on the full measure subset.
\end{itemize}

Moreover in classical mechanics the most popular situation is when
$M^{2n}$ and $H$ are (real-)analytic and

\begin{itemize}
\item
first integrals $I_1,\dots,I_n$ are analytic.
\end{itemize}

In this case it is said that the flow is analytically integrable.

Another reasonable treatment of what means ``almost every functional
independence of integrals of motion''
for a compact phase space is as follows

\begin{itemize}
\item
there is a finite smooth (or even analytic)
simplicial decomposition of the phase space $M^{2n}$ such that
a singular locus $\Sigma$ forms a subcomplex of this decomposition
and the complement to it is cutted by another subcomplex of positive
codimension to a union of finitely many sets $U_\alpha$
which are foliated by invariant tori over their
images under the momentum map.
\end{itemize}

We considered such a notion in \cite{T1} and called
it a geometric simplicity.

It is reasonable to say that integrals of motion are tame if they lead
to a geometrically simple behaviour of the flow.

Some important examples of Hamiltonian flows do not
have a compact phase space.
This is, for instance, the geodesic flow of a Riemannian manifold $M^n$
which is a Hamiltonian flow on a cotangent bundle
$T^\ast M^n$ to this manifold.
The symplectic structure on $T^\ast M^n$ is given by a form
$$
\omega = \sum_{j=1}^n dx^j \land dp_j
$$
where $p_j = g_{jk} \dot{x}^k$ and $g_{jk}dx^j dx^k$
is the Riemannian metric.
The Hamiltonian of the geodesic flow is homogeneous in momenta:
$$
H(x,p) = g^{jk}(x) p_j p_k
$$
where $g^{jk}g_{kl} = \delta^j_l$ and therefore
the restrictions of the geodesic
flow onto different non-zero level surfaces of $H$ are smoothly trajectory
equivalent, moreover they are related by reparametrization of trajectories
by a constant. Therefore it is reasonable to call the geodesic flow on $M^n$
integrable if it satisfies a weaker condition (see \cite{T1})
which is

\begin{itemize}
\item
there are $(n-1)$ additional integrals of motion $I_1,\dots,I_{n-1}$
which are in involution and almost everywhere
independent on the unit covector
bundle $SM^n = \{H=1\}$.
\end{itemize}

It appears that an analytic integrability looks the strongest assumption which
implies, in particular, geometric simplicity as it was shown in \cite{T1}.

Some recent examples of integrable geodesic flows of analytic metrics
show that even using for integration of geodesic flows
such mildly non-analytic $C^\infty$ functions as, for
instance,
$$
f(x,p) = \exp(- Q(p)^{-2}) \sin(\mu \log|p_x - \tau p_y|)
$$
where $\mu, \tau$ are constants and $Q(p)$ is
a polynomial in momenta $p_x,p_y$
divided by $(p_x - \tau p_y)$, leads to
geometrically non-simple flows \cite{BT,Butler}.

Therefore

{\it for considering topological properties of integrable flows by means
of topology of finite CW-complexes or tame subsets in $\R^n$ we have to
restrict ourselves to geometrically simple flows or to tame integrals of
motion.}

In the next chapter we expose some background from mathematical logics which
leads to the most clarified analytic approach to understanding what it
means that an integral of motion is tame.

Before we do that we would like to notice that
an analogue of the Morse theory
for integrable systems on four-dimensional symplectic manifolds developed by
Fomenko and his collaborators also needs some analytic condition
fulfilled. It reads that an additional (to the Hamiltonian)
integral of motion has to be of the Bott type, i.e. its restrictions onto
normal planes to critical level surfaces would be locally a Morse function.
This is another kind of tameness condition which was generalized in
\cite{MF} to a geometric condition under which this theory works.

\section{O-minimal structures and analytic-geomet\-ric
ca\-te\-go\-ries}

\subsection{O-minimal structures}

By definition, a
family ${\cal S}$
of subsets of the Euclidean spaces $\R^n$ is called an o-minimal structure
(on $(\R,+,\cdot,<)$) if being graded by the
dimensions of ambient Euclidean spaces:
$$
{\cal S} = \cup_{n \geq 1} {\cal S}_n,
$$
it meets the following conditions:

1) ${\cal S}_n$ is a Boolean algebra of some subsets of $\R^n$
with a standard union operation (in particular, this means that this algebra
is closed with respect to complements and finite unions and intersections);

2) if $X \in {\cal S}_n$ and $Y \in {\cal S}_k$,
then $X \times Y \in {\cal S}_{n+k}$;

3) let $\pi: \R^{n+1} \to \R^n$ be a projection $(x^1,\dots,x^n,x^{n+1}) \to
(x^1,\dots,x^n)$, then $X \in {\cal S}_{n+1}$
implies $\pi(X) \in {\cal S}_n$;

4) ${\cal S}_n$ contains all algebraic sets in $\R^n$, i.e. if
$P(x^1,\dots,x^n)$ is a polynomial then its zero set $\{P =0\}$ belongs to
${\cal S}$;

5) ${\cal S}_1$ consists exactly
of all finite unions of open intervals and points.

Let us present the main known examples of o-minimal structures:

\begin{itemize}

\item
$\R_{\mathrm{alg}}$: {\bf semialgebraic sets}. It consists exactly
of all semialgebraic sets, i.e. determined by finitely many equations
$P_1 = \dots = P_k = 0$, and inequalities
$Q_1 >0$, $\dots$, $Q_l > 0$ with $P_1,\dots,P_k,Q_1,\dots,Q_l$
polynomials. Such sets form an o-minimal structure by the
Tarski--Seidenberg theorem.

\item
$\R_{\mathrm{an}}$:
{\bf finite subanalytic sets}. It consists of intersections
of subanalytic sets in $\R^n$ with cubes $[-D,D]^n$ and their projections.
Notice that the family formed only by intersections of subanalytic
sets with cubes is not closed under projections. It is the theorem of
Gabrielov \cite{Gabrielov} which implies that $\R_{\mathrm{an}}$ is an
o-minimal structure.

\item
$\R_{\exp}$. Let us correspond to any polynomial $P(x^1,\dots,
x^{2k})$ an exponential polynomial $Q(x^1,\dots,x^k) =
P(x^1,\dots,x^k,e^{x^1},\dots,e^{x^k})$ and
denote by $\R_{\exp}$ the family of
all sets generated by the zero sets of such exponential polynomials under
projections $\R^{n+1} \to \R^n$. Wilkie proved that this family is
closed under complements and forms an o-minimal structure \cite{W1}.

\item
$\R_f$, where $f$ is a Pfaffian function. This is a generalization of
$\R_{\exp}$. It is said that a chain of $C^1$-functions
$f_1,\dots,f_k: \R^n \to \R$ is a Pfaffian chain if for each $j=1,\dots,k$
the first derivatives of $f_j$ with respect to $x^1,\dots,x^n$
are polynomials in $x^1,\dots,x^n,f_1,\dots,f_j$. In this event the function
$f=f_k$ is called a Pfaffian function.
The similarity between the zero sets of
Pfaffian functions and algebraic sets was first pointed out by Khovanskii
\cite{Kh}. It was proved by Wilkie that if we replace $\exp$ by a Pfaffian
function $f$ in the definition of $\R_{\exp}$ and close this family with
respect to Boolean operations, projections and products then we get an
o-minimal structure which is denoted by $\R_f$ \cite{W2}.

\item
$\R_{\mathrm{an,exp}}$.
It was proved by van den Dries, Macintyre, and Marker \cite{DMM} by using
the results of Wilkie \cite{W1}
that sets from $\R_{\mathrm{an}}$ and
$\R_{\exp}$ generate by Boolean operations and projections
an o-minimal structure which is denoted by $\R_{\mathrm{an,exp}}$.
\end{itemize}

These o-minimal structures are related by the following evident inclusions:
$$
\R_{\mathrm{alg}} \subset \R_{\mathrm{an}} \subset \R_{\mathrm{an,exp}}.
$$

Given an o-minimal structure ${\cal S}$, we say that

\begin{itemize}

\item
a subset $X \in \R^n$
is definable if $X \in {\cal S}_n$
\footnote{This terminology originates in mathematical logics and reflects the
fact that definiable sets are exactly sets which are defined by the first
order logics formulas involving the summation ``$+$", the 
multiplication ``$\cdot$'' and the linear ordering ``$<$'' plus some additional
functions which lead to extensions of 
the smallest o-minimal structure on $(\R,+,\cdot,<)$, i.e.
the subalgebraic sets. These are analytic functions restricted to cubes
$[-1,1]^n$ for $\R_{\mathrm{an}}$, the exponent function $\exp$ for 
$\R_{\mathrm{an}}$ and etc. If we drop the multiplication from the signature 
of our language (in the sense of mathematical logics) 
we have to replace the 4th condition by another which reads that the graphs 
of some functions coming into the signature are definable. For instance,
the smallest o-minimal structure which includes $+$, $<$ and
the multiplications $r$ by all real numbers $r \in \R$ is formed by all
semilinear sets.};

\item
a mapping $f: X \to \R^k$ with $X \subset \R^n$ is definable if
its graph $\{(x,f(x))\}$ is a definable set,
i.e. belongs to ${\cal S}_{n+k}$;

\item
a set $X$ is ${\cal S}$-triangulable if
$X \in {\cal S}$ and there is a definable mapping
$f: X \to \R^n$ which maps $X$
homeomorphically onto a union of open simplices of finite simplicial
complex $K \subset \R^n$. In this event we say that
$f$ defines an ${\cal S}$-triangulation of $X$.
\end{itemize}

By this definition, any ${\cal S}$-triangulable set
is definable. The converse is also true:

\begin{theorem}[Triangulation Theorem]
Every definable set $X \in {\cal S}$ is ${\cal S}$-tri\-an\-gu\-lable.
\end{theorem}

This theorem is proved by a general method for all o-minimal structures
\cite{Dries}. The proof follows by induction on the dimension of a
definable set
and starts with an evident statement that all sets from ${\cal S}_1$ are
${\cal S}$-triangulable. We sketched such a proof for sets from
$\R_{\mathrm{an}}$ in \cite{T1}.

Let us also notice that it follows from the definition of
an o-minimal structure that
images and preimages of definable sets under definable proper mappings
are definable. Here we recall that a mapping is called proper if
preimages of compact sets are compact.

\subsection{Geometric and analytic-geometric categories}

For using the theory of o-minimal structures in topology and geometry
it needs to develop its analog for subsets in manifolds and it was done
in \cite{DM}.

Given an o-minimal structure ${\cal S}$, we distinguish a class of
${\cal S}$-manifolds. We say that a smooth manifold $M^n$ is a
${\cal S}$-manifold if it
admits a finite ${\cal S}$-atlas $\{U_\alpha\}$, i.e.
such an atlas that

\begin{itemize}
\item
every coordinate
mapping $\varphi_\alpha: U_\alpha \to V_\alpha \subset \R^n$
homeomorphically maps a chart $U_\alpha$ onto a definable set $V_\alpha$;

\item
for any intersection $U_\alpha \cap U_\beta$ the transition mapping
$$
\varphi_\alpha \varphi_\beta^{-1}:
\varphi_\beta(U_\alpha \cap U_\beta) \to
\varphi_\alpha(U_\alpha \cap U_\beta)
$$
is definable.
\end{itemize}

Now we say that

\begin{itemize}
\item
a subset $X \in M^n$ is definable if for any chart $U_\alpha$
the set $\varphi_\alpha (X \cap U_\alpha) \subset \R^n$ is definable;

\item
for definable sets $X \in M^n$ and $Y \in N^k$ a mapping $f: X \to Y$
is definable if it is continuous and its graph $\{(x,f(x))\}$ is definable
in $M^n \times N^k$.
\end{itemize}

Notice that these definitions are independent on choosing ${\cal S}$-atlases
for $M^n$ and $N^k$.

Thereby an o-minimal structure ${\cal S}$
defines a ``geometric category'' of
${\cal S}$-manifolds and their definable subsets as objects and
definable maps between them as morphisms. For instance, if ${\cal S} =
\R_{\mathrm{alg}}$ we have a category of
so-called Nash manifolds \cite{Shiota}.

To an o-minimal structure ${\cal S}$ which contains $\R_{\mathrm{an}}$ there
corresponds an ``ana\-ly\-tic-geometric'' category which is defined as
follows

\begin{itemize}
\item
given an analytic manifold $M^n$ a subset $X \subset M^n$ is called
definable if for any point $x \in X$ there are an open neighborhood $U$
of this point and an analytic isomorphism
$\varphi: U \to V \in \R^n$ such that
$\varphi(X \cap U)$ is definable;

\item
a mapping $f: X \to Y$ of two definable sets
$X \subset M^n$ and $Y \subset N^k$
is called definable if its graph is definable in $M^n \times N^k$.
\end{itemize}

Notice that by this definition an image of a definable set
under a proper analytic mapping is a definable set.

The analytic-geometric category corresponding to ${\cal S}$ has
definable sets as objects and definable mappings as morphisms.

Introduction of these categories allows us to use machinery developed for
definable sets in $\R^n$, for instance Triangulation Theorem and many
other facts (see \cite{DM,Dries}), for subsets in manifolds.

\section{On obstructions to integrability}

The first obstruction to integrability of geodesic flows was found by
Kozlov who proved that if the geodesic flow on a compact oriented
analytic Riemannian two-manifold admits an additional analytic
first integral then the surface is diffeomorphic to a sphere or to a torus
\cite{Kozlov}.

The analycity condition was strongly used in the proof
and it is still unknown
is the theorem is valid under only $C^\infty$ assumptions for a manifold and
first integrals. Using the theory of modular forms Kolokoltsov extended
the Kozlov theorem assuming
that an additional first integral is $C^\infty$ but
polynomial in momenta \cite{Kol}.

High-dimensional obstructions were obtained by us in \cite{T1} in two steps:

1) there were found some obstructions to a nice geometric behaviour
of the geodesic flow on a manifold, i.e. obstructions to its geometric
simplicity,

2) some analytic properties of first integrals
which imply geometric simplicity were established.

We shall remark here that the condition of geometric simplicity
can be weakened
and by using the language of o-minimal structures the analytic condition
can be clarified and slightly strengthened.

For realizing the first step we proved the following

\begin{theorem}
If the geodesic flow on a compact analytic manifold $M^n$ is
geometrically simple then there is an invariant torus $T^n \subset
SM^n$ such that the natural projection $\pi: SM^n \to M^n$ induced
a homomorphism
$$
\pi_\ast: \pi_1(T^n) \to \pi_1(M^n)
$$
those image is a subgroup of finite index in $\pi_1(M^n)$.
\end{theorem}

It implies

\begin{corollary}
\label{topology}
If a geodesic flow on a compact manifold $M^n$ is
geometrically simple then

1) the fundamental group $\pi_1(M^n)$ of $M^n$ is almost commutative,
i.e., contains a commutative subgroup of finite index;

2) the real cohomology ring $H^\ast(M;\R)$ of $M^n$ contains a subring $A$
which is isomorphic to the real cohomology ring $H^\ast(T^k;\R)$ of the
$k$-dimensional torus where $k$ is the first Betti number of $M^n$:
$b_1 = \dim H^1(M^n;\R) = k$;

3) moreover if the first Betti number of $M^n$ equals its dimension: $b_1
=n$ then the ring $H^\ast(M^n;\R)$ is isomorphic to $H^\ast(T^n;\R)$.
\end{corollary}

To explain these results we recall that we say that a geodesic flow on $M^n$
is geometrically simple if the unit cotangent bundle $SM^n$
admits a decomposition
$$
SM^n = \Gamma \cup \left(\cup_{\alpha=1}^d U_\alpha \right)
$$
such that

\begin{itemize}
\item
this decomposition is invariant under the flow;

\item
the set $\Gamma$ is closed and the complement to it is everywhere dense;

\item
for any point $q \in SM^n$ and every its neighborhood $V$ there is
another neighborhood $W$ of $q$
such that $W \subset V$ and $W \cap (M^n \setminus \Gamma)$
has finitely many connected components;

\item
any component $U_\alpha$ is diffeomorphic to a product of an $n$-dimensional
torus and an $(n-1)$-dimensional disc.
\end{itemize}

In fact we proved that omitting the fourth condition there is
a component $U_\alpha$ such that the image of its fundamental group under
the projection homomorphism $\pi_\ast(\pi_1(U_\alpha))$
has a finite index in $\pi_1(M^n)$.

Moreover in this formulation the proof of the theorem works for
more general case when the flow is locally simple, i.e.
there is a point $x \in M^n$ and its neighborhood $U$ such that

\begin{itemize}
\item
the universal covering $\widehat{M}^n \to M^n$ is trivial over $U$;

\item
the preimage of $U$ under the projection $\pi: SM^n \to M^n$ admits a
decomposition
$$
\pi^{-1}(U) = \widetilde{\Gamma} \cup
\left( \cup_{\alpha =1}^d \widetilde{U}_\alpha \right)
$$
where $\widetilde{\Gamma}$ is closed and the complement to
it is dense, and each
component $\widetilde{U}_\alpha$ is an intersection of $\pi^{-1}(U)$ with an
invariant open set $U_\alpha)$.
\end{itemize}

The second step was realized in \cite{T1} by the following

\begin{theorem}
\label{analytic}
If a geodesic flow on a compact manifold is analytically integrable then
it is geometrically simple.
\end{theorem}

In proving this theorem
the basic point is to show  that given analytic first integrals
$I_1,\dots,I_{n-1}$ (here we assume
that $I_n$ is the Hamiltonian of the flow,
$I_n = g^{ij}(x)p_i p_j$) the set $C$
of the critical values of the momentum map restricted onto $SM^n$
$$
\Phi: q \to (I_1(q),\dots,I_{n-1}(q)) \in \R^{n-1}
$$
and its preimage in $SM^n$ are analytically-triangulable.

In the modern terminology of {\S} 3,
the proof of that consist in a remark that these sets
$C$ and $\Phi^{-1}(C)$ are definable in
$\R_{\mathrm{an}}$-analytic-geometric category and therefore are
$\R_{\mathrm{an}}$-triangulable. We proved their
analytic-trian\-gu\-la\-bi\-lity
directly by using the Gabrielov theorem \cite{Gabrielov}.
We already mentioned that the proof of Triangulation Theorem for general
o-minimal structures follows the same scheme as we used
which probably originates in Hironaka's proof of
Triangulation Theorem for semialgebraic sets.

After proving that $C$ and $\Phi^{-1}(C)$ are analytically-triangulable
we completed the set $C$
of $\Phi$ by adding some additional analytic $(n-2)$-dimensional simplices
to simplical subcomplex $K$ those complement in $\Phi(SM^n)$ is a union of
finitely many discs $V_\alpha$ and denote $\Phi^{-1}(V_\alpha)$
by $U_\alpha$ thus proving a geometric simplicity.

Now by using the general form triangulation Theorem we can
generalize this theorem as follows:

\begin{theorem}
\label{structure}
Let ${\cal S}$ be an o-minimal structure.
Let $M^n$ is a compact
Riemannian ${\cal S}$-manifold and assume that the geodesic
flow on $M^n$ is integrable in terms of ${\cal S}$-definable first integrals.
Then this geodesic flow is geometrically simple.
\end{theorem}

For ${\cal S} = \R_{\mathrm{an}}$
this theorem reduces to Theorem \ref{analytic}.

It was first shown by Butler that assuming only integrability in terms of
$C^\infty$ first integrals we can not conclude
that the fundamental group of the
manifold is almost commutative. He did that by constructing a $C^\infty$
integrable geodesic flow on a three-dimensional nilmanifold \cite{Butler}.
Later Bolsinov and the author even managed to construct a $C^\infty$
integrable geodesic flow on a solvmanifold those fundamental group has
an exponential growth \cite{BT}.

But as Theorem \ref{structure} shows we can derive topological
conclusions of Corollary \ref{topology} by assuming
that the flow is integrated
in terms of $C^\infty$ first integrals which are definable in some
analytic-geometric category. In this event the category
corresponding to $\R_{\mathrm{an}}$ is the smallest possible category.


\begin{thebibliography}{MMM}

\bibitem{Arnold}
Arnold,V.I.
Mathematical methods of classical mechanics.
Springer-Verlag, New York, 1989.

\bibitem{BT}
Bolsinov, A.V., Taimanov, I.A.
Integrable geodesic flows with positive topological entropy.
Inventiones Math. {\bf 140} (2000), 639--650.

\bibitem{Butler}
Butler, L.T.
New examples of integrable geodesic flows.
Asian J. Math. {\bf 4} (2000), 515--526.

\bibitem{DM}
van den Dries, L., Miller, C.
Geometric categories and o-minimal structures.
Duke Math. J. {\bf 84} (1996), 497--540.

\bibitem{Dries}
van den Dries, L.
Tame topology and o-minimal structures.
London Math. Society Lecture Notes {\bf 248},
Cambridge University Press, Cambridge, 1998.

\bibitem{DMM}
van den Dries, L., Macintyre, A., Marker, D.
The elementary theory of restricted analytic fields with exponentiation.
Annals of Math. {\bf 140} (1994), 183--205.

\bibitem{Gabrielov}
Gabrielov, A.M.
Projections of semianalytic sets.
Functional Anal. Appl. {\bf 2} (1968), 282--291.

\bibitem{Kh}
Khovanskii, A.G.
Fewnomials.
Translations of Math. Monographs, Amer. Math. Soc. {\bf 88},
Providence, 1991.

\bibitem{Kol}
Kolokoltsov, V.N.
Geodesic flows on two-dimensional manifolds
with an additional first integral
that is polynomial with respect to velocities.
Math. USSR-Izv. {\bf 46} (1983), 291--306.

\bibitem{Kozlov}
Kozlov, V.V.
Topological obstacles to the integrability of natural mechanical systems.
Soviet Math. Dokl. {\bf 20} (1979), 1413--1415.

\bibitem{MF}
Matveev, S.V., Fomenko, A.T.
A Morse-type theory for integrable Hamiltonian systems with tame integrals.
Math. Notes {\bf 43}, no. 5-6 (1988), 382--386.

\bibitem{Shiota}
Shiota M. Nash Manifolds. Lecture Notes in Math. {\bf 1269},
Springer-Verlag, New-York, 1987.

\bibitem{T1}
Taimanov, I.A.
Topological obstructions to the integrability of geodesic flows on
nonsimply connected manifolds.
Math. USSR-Izv. {\bf 30} (1988), 403--409.

\bibitem{T2}
Taimanov, I.A.
Topology of Riemannian manifolds with integrable geodesic flows.
Proc. Steklov Inst. Math. {\bf 205} (1995), 139--150.

\bibitem{W1}
Wilkie, A.J.
Model completeness results for expansions of the ordered field of
real numbers by restricted Pfaffian functions and the exponential function.
J. of the Amer. Math. Soc. {\bf 9} (1996), 1051--1094.

\bibitem{W2}
Wilkie, A.J.
A theorem of the complement and some new o-minimal structures,
Selecta Math. (N.S.) {\bf 5} (1999), 397--421.


\end{thebibliography}
\end{document}